\documentclass[10pt]{article}
\usepackage{amsmath,amsthm,amscd,amssymb}
\newtheorem{df}{Definition}[section]
\newtheorem{thm}{Theorem}[section]
\newtheorem{prop}{Proposition}[section]
\newtheorem{lm}{Lemma}[section]

\newtheorem{fact}{Fact}[section]
\newtheorem{cor}{Corollary}[section]

\title{A special value of Ruelle L-function and the theorem of Cheeger
and M\"{u}ller}
\author{Ken-ichi SUGIYAMA\thanks{Cooresponding address : Department of
Mathematics and Informatics, Faculty of Science, Chiba University, 1-33
Yayoi-cho Inage-ku, Chiba 263-8522, Japan. e-mail address : sugiyama@math.s.chiba-u.ac.jp}\\
Department of Mathematics and Informatics, Faculty of Science\\
Chiba University, Japan.}

\begin{document}
\maketitle
\begin{abstract}
We will show a theorem of a type of Cheeger and M\"{u}ller for a noncompact
 complete hyperbolic threefold of finite volume. As an application we
 will compute a special value of Ruelle L-function at the origin for a
 unitary local system which is cuspidal. 
\par\vspace{5pt}
\end{abstract}
\section{Introduction}
A special value of an L-function associated to a representation of the absolute Galois group of a
number field reflects an arithmetic or geometric property of the base
field or an object from which the representation arises. For example the
class number formula says that Dedekind zeta function $\zeta_F(s)$ of a
number field $F$, which is an L-function associated to the trivial representation, has a simple pole at $s=1$ and the residue is expressed
in terms of arithmetic invariant of $F$, e.g. a class number, a
fundamental regulator and a number of roots of unity in $F$ and so
on. Birch and Swinnerton-Dyer conjecture for an elliptic curve defined
over ${\mathbb Q}$ predicts that an L-function of associated $l$-adic
representation should have zero at $s=1$ whose order is equal to the rank
of Mordell-Weil group $E({\mathbb Q})$. Moreover it says that the
leading coefficient of Taylor expansion at $s=1$ should be written by arithmetic/geometric invariants
of $E$, e.g. an order of Shafarevich-Tate group, an elliptic regulator
and an order of the torsion subgroup of $E({\mathbb Q})$ and so on. \\

In the present paper we will discuss an analog of such formulas for a
unitary representation $\rho$ whose degree $r$ of the fundamental group $\pi_1(X)$ of a complete
hyperbolic threefold $X$ of finite volume. Such a
representation associates a unitary local system on $X$, which will be
denoted by the same character. We assume that a restriction of $\rho$ to
a fundamental group at every cusp does not fix any vector other than
$0$. (If this is satisfied we call $\rho$ {\it cuspidal}.)\\

 In order to
explain Ruelle L-function we prepare some terminologies. Since $X$ is
hyperbolic $\pi_1(X)$ may be identified with a discrete subgroup of
${\rm PSL}_2({\mathbb C})$ and there is the natural bijection between a
set of hyperbolic conjugacy classes of $\pi_1(X)$ and a set of closed
geodesics of $X$. Using this the length $l(\gamma)$ of a hyperbolic
conjugacy class $\gamma$ is defined to be one of the corresponding
geodesic. A closed geodesic will be referred as {\it prime} if it is not
a positive multiple of an another one. Using the bijection we define a subset
$\Gamma_{prim}$ of hyperbolic conjugacy classes which
consists of elements corresponding to prime closed geodesics.  
Now Ruelle L-function is defined to be
\[
 R_{X}(z,\rho)=\prod_{\gamma\in \Gamma_{prim}}\det[1-\rho(\gamma)e^{-zl(\gamma)}]^{-1}.
\]
It absolutely converges if ${\rm Re}s>2$. J.Park has shown that it is
meromorphically continued to the whole plane and that it has zero at the
origin of order $2h^1(X, \rho)$, where $h^p(X, \rho)$ is the dimension of
$H^p(X, \rho)$(\cite{Park2007}). We will show the following theorem.
\begin{thm}
\[
 \lim_{z\to 0}z^{-2h^{1}(X,\rho)}R_{X}(z,\rho)=(\tau^{*}(X,\rho)\cdot{\rm Per}(X))^2.
\]
\end{thm}
In the theorem $\tau^{*}(X,\rho)$ is a modified Franz-Reidemeister
torsion (see \S4) and ${\rm Per}(X)$ is a period of $X$ (see \S5). The
former is a combinatric invariant and the latter is an analytic
one. Notice that {\bf Theorem 1.1} may be compared to Birch and
Swinnerton-Dyer conjecture.
\begin{cor} Suppose that $h^{1}(X,\rho)$ vanishes. Then
\[
 R_{X}(0,\rho)=\tau(X,\rho)^{2},
\]
where $\tau(X,\rho)$ is the usual Franz-Reidemeister
torsion.
\end{cor}
Here is an example so that $R_{X}(0,\rho)$ is computed explicitly. Let
$K$ be a hyperbolic knot in $S^3$. Thus its complement $X_K$ admits
hyperbolic structure of finite volume. Let $\xi$ be a complex number of
modulus one. Since $H_1(X_K,{\mathbb Z})$ is isomorphic to an infinite
cyclic group, sending a generator to $\xi$, we have a map
\[
 H_1(X_K,{\mathbb Z}) \to {\rm U}(1),
\]
and composing with Hurewicz map it induces a unitary character
\[
 \pi_1(X_K) \stackrel{\rho_{\xi}}\to {\rm U}(1).
\]
It is easy to see that if $\xi\neq 1$, $\rho_{\xi}$ is
cuspidal. Moreover we can show that the fact $A_K(\xi)\neq 1$ is equivalent
to $h^{1}(X,\rho)=0$, where $A_K(t)$ is Alexander polynomial.
\begin{cor} Let us choose $\xi$ so that
both $\xi-1$ and $A_{K}(\xi)$ do not vanish. Then
\[
 R_{X_K}(0,\rho)=\left|\frac{A_{K}(\xi)}{1-\xi}\right|^2.
\]
\end{cor}
A proof of {\bf Theorem 1.1} is based on a result of J. Park(\cite{Park2007}) which
asserts the leading coeffcient of Taylor expansion of $R_X(z,\rho)$ at
the origin is $\exp(-\zeta^{\prime}_X(0,\rho))$, where
$\zeta_X(z,\rho)$ is the spectral zeta function (see \S4). Thus {\bf
Theorem 1.1} is reduced to show an equation:
\[
\exp(-\zeta^{\prime}_X(0,\rho))= (\tau^{*}(X,\rho)\cdot {\rm Per}(X))^2,
\]
or equivalently to show the following theorem of Cheeger-M\"{u}ller
type.
\begin{thm}
$||\cdot||_{FR}$ and $||\cdot||_{RS}$ coincide.
\end{thm}

First the theorem has been independently proved by
Cheeger(\cite{Cheeger}) and M\"{u}ller(\cite{Muller}) for a closed
manifold, which are solutions of Ray-Singer conjecture. For a compact manifold with boundaries, if a metric is a
product near boundaries, it has been independently observed by
Lott-Rothenberg(\cite{Lott-Rothenberg}), L\"{uck}(\cite{Luck}) and
Vishik(\cite{Vishik}) that $||\cdot||_{FR}$ and
$||\cdot||_{RS}$ differ by Euler characteristic of the restriction of
$\rho$ to boundaries. Moreover Dai and Fang (\cite{Dai-Fang}) have
computed their difference when a metric is not a product near
ends. In our case, cutting by holospheres, $X$ may be considered
as a limit of compact Riemannian manifolds with torus boundaries whose metric is
not a product near ends. Using results of Dai and Fang we will estimate difference
between $||\cdot||_{FR}$ and $||\cdot||_{RS}$ and will show that their
limit coincide.\\

{\bf Acknowledgment.} The author express heartly gratitude to Professor
J. Park who kindly show him a preprint \cite{Park2007}, which is
indispensable to finish this work.

\section{Spectrum of Laplacian near cusps}
Let $X$ be a complete hyperbolic threefold of finite volume with cusps
$\{\infty_{\nu}\}_{1\leq \nu\leq h}$. Thus it is a quotient of
Poincar\'{e} upper half space ${\mathbb H}^3=\{(x,y,r)\in{\mathbb
R}^3\,|\,r>0\}$ equipped with a metric
\[
 g=\frac{dx^2+dy^2+dr^2}{r^2}
\]
of contant curvature $-1$ by a discrete subgroup $\Gamma$ of ${\rm
PSL}_2({\mathbb C})$. A cusp $\infty_{\nu}$ associates the unipotent
radical $N_{\nu}$ of a Borel subgroup $B_\nu$ of ${\rm PSL}_2({\mathbb
C})$. Without loss of generality we may assume $B_1$ is the standard
Borel subgroup which consists of upper triangular matrices and thus
\[
 N_1=\{
\left(
\begin{array}{cc}
1&z\\
0&1\end{array}
\right)\,|\,z\in{\mathbb C}\}.
\]
For $1\leq \nu \leq h$ there is $g_{\nu}\in{\rm PSL}_2({\mathbb C})$
such that
\[
 B_{\nu}=g_{\nu}B_1g_{\nu}^{-1} \quad \mbox{and} \quad N_{\nu}=g_{\nu}N_1g_{\nu}^{-1}.
\]
Here we take $g_1$ to be the identity matrix.
Now {\it the fundamental group at} $\infty_{\nu}$ is defined to be
\[
 \Gamma_{\nu}=\Gamma\cap N_{\nu},
\]
which is isomorphic to ${\mathbb Z}\oplus{\mathbb Z}$.\\

For a positive $a$ we put
\[
 {\mathbb H}^{3}_{a,\infty}=\{(x,y,r)\in{\mathbb H}^3\,|\, r \leq e^{a}\}
\]
and
\[
 {\mathbb H}^{3}_{a}=\cap_{\nu=1}^{h}g_{\nu}{\mathbb H}^{3}_{a,\infty}.
\]
Let $X_{a}$ be the image of ${\mathbb H}^{3}_{a}$ by the natural projection and $Y_{a}$ 
the closure of its complement. If $a$ is sufficiently large $Y_{a}$ is a
disjoint union of $Y_{a,\nu}\,(1\leq \nu \leq h)$ and an each of them is topologically a product of a flat 2-torus
$T_{\nu}=N_{\nu}\slash \Gamma_{\nu}\simeq {\mathbb C}\slash \Gamma_{\nu}$ and
an interval $[e^a,\infty)$. Moreover by a change of variables
\[
 r=e^{u},
\]
$Y_{a,\nu}$ becomes a warped product $[a,\infty)\times T_{\nu}$ equipped
with a metric
\[
 g=du^2+e^{-2u}(dx^2+dy^2).
\]
In particular the boundary of $Y_{a,\nu}$ is $T_{\nu}$ but with a metric $e^{-2a}(dx^2+dy^2)$.
In the following computations will be carried out by a coordinate
$(x,y,u)$.\\

Let $\rho$ be a unitary representation of $\Gamma$ of rank $r$. It
yields a unitary local system on $X$ which will be denoted by the same
character. Since $\Gamma_{\nu}$ is abelian the restriction
$\rho|_{\Gamma_{\nu}}$ is decomposed into a direct sum of unitary
characters:
\begin{equation}
\rho|_{\Gamma_{\nu}}=\oplus_{i=1}^{r}\chi_{\nu,i}.
\end{equation}
Throughuot the paper we will always assume that $\rho$ is {\it cuspidal} i.e. none of
$\{\chi_{\nu,i}\}_{\nu,i}$ is trivial. A vector bundle of $p$-forms on $X$
twisted by $\rho$ will be denoted by $\Omega^{p}_{X}(\rho)$. More
generally for a submanifold $M$ let $\Omega^{p}_{M}(\rho)$ be a vector
bundle of $p$-forms on $M$ twisted by $\rho$. Let $\varphi$ be a smooth
section of $\Omega^{p}_{X}(\rho)$ on $Y_{a,\nu}$. By decomposition (1)
we have
\[
 \varphi=\sum_{i=1}^{r}\varphi_i,\quad \varphi_i=\sum_{|\alpha|=p}\varphi_{i,\alpha}dx^{\alpha}\in
 C^{\infty}(Y_{a,\nu}, \Omega_{X}^{p}(\chi_{\nu,i})).
\]
Here we have used a convention:
\[
 x_0=u,\quad x_1=x, \quad \mbox{and} \quad x_2=y.
\]
\begin{lm}
$\varphi$ is cuspidal, i.e. for any $\nu$, $i$ and $\alpha$ 
\[
 \int_{T_{\nu}}\varphi_{i,\alpha}dxdy=0.
\]
\end{lm}
{\bf Proof.} Let us choose $\gamma\in \Gamma_{\nu}$ so that
\[
 \chi_{\nu,i}(\gamma)\neq 1.
\]
By definition we have
\[
 \gamma^{*}\varphi_{i,\alpha}=\chi_{\nu,i}(\gamma)\varphi_{i,\alpha}.
\]
and the desired result will follow from
\[
 \int_{T_{\nu}}\varphi_{i,\alpha}dxdy=\int_{T_{\nu}}\gamma^{*}\varphi_{i,\alpha}dxdy=\chi_{\nu,i}(\gamma)\int_{T_{\nu}}\varphi_{i,\alpha}dxdy.
\]
\begin{flushright}
$\Box$
\end{flushright}
We will consider an eigenvalue problem of Hodge Laplacian $\Delta^p$ on
spaces of square integrable twisted $p$-forms 
$L^2(X_a,\Omega^p_{X}(\rho))$ or
$L^2(Y_{a,\nu},\Omega^p_{X}(\rho))$ under {\it an absolute},
{\it a relative} or {\it Dirichlet} boundary condition, which we will
now recall. The restriction $\Omega_{X}^p(\rho)$ to the boundary
$T_{\nu}$ of $Y_{a,\nu}$ is decomposed as
\[
 \Omega_{X}^p(\rho)|_{T_{\nu}}=\Omega^{p}_{T_{\nu}}(\rho)\oplus du\wedge \Omega^{p-1}_{T_{\nu}}(\rho).
\]
According to this a section $\omega$ of a restriction of
$\Omega^p_{X}(\rho)$ to $T_{\nu}$ is written to be
\[
 \omega=\omega_{tan}+ \omega_{norm},
\]
where $\omega_{tan}$ (resp. $\omega_{norm}$) is a section of
$\Omega^{p}_{T_{\nu}}(\rho)$ (resp. $du\wedge\Omega^{p-1}_{T_{\nu}}(\rho)$).
\begin{df} We call $\omega$ satisfies {\rm an absolute boundary condition} if
 both $\omega_{norm}$ and $(d\omega)_{norm}$ vanish on every
 connected component $T_{\nu}$ of the boundary. If the Hodge dual $*\omega$ satisfies an absolute boundary condition $\omega$
 will be referred as it satisfies {\rm a relative boundary
 condition}. More strongly if the restrictions of both $\omega$ and $d\omega$ to
 $T_{\nu}$ vanish for every $\nu$, we call it satisfies {\rm Dirichlet
 boundary condition}.
\end{df}
Notice that $*$ interchanges the first two conditions and preserves the
last one. Since $\rho$ is unitary the local system possesses a fiberwise
hermitian inner product ${\rm Tr}_{\rho}$. For $\omega,\,\eta\in
\Omega^p_{X}(\rho)$ we put
\[
 (\omega,\eta)=\frac{{\rm Tr}_{\rho}(\omega\wedge *\eta)}{dv_g},
\]
which becomes a hermitian inner product on $\Omega_{X}^p(\rho)$. Here
$dv_g$ is the volume form of $g$, which is equal to $e^{-2u}dx\wedge
dy\wedge du$.
Let $M$ be $X_{a}$ or $Y_{a,\nu}$. If both of $\omega$ and $\eta$
satisfies one of boundary conditions we have by Stokes theorem
\[
 \int_{M}(\Delta^p\omega,\eta)dv_g=\int_M(\nabla\omega,\nabla\eta)dv_g=\int_{M}(\omega,\Delta^p\eta)dv_g,
\]
where $\nabla$ is the covariant derivative. Therefore $\Delta^p$ has a
selfadjoint extension $\Delta^p_{abs}$, $\Delta^p_{rel}$ or
$\Delta^p_{dir}$ according to a boundary condition which is absolute,
relative or Dirichlet, respectively. If $\sharp$ is $abs$ (resp.
$rel$ or $dir$) {\it its dual} $\hat{\sharp}$ is defined to be $rel$ (resp.
$abs$ or $dir$). Since by Hodge symmetry a Hilbert module
$\{L^2(M,\Omega^p_X({\rho})),\,\Delta^p_{\sharp}\}$ is isomorphic to
$\{L^2(M,\Omega^{3-p}_X({\rho})),\,\Delta^{3-p}_{\hat{\sharp}}\}$ we
will only consider the case of $p=0$ or $1$. For a later purpose we will
introduce one more boundary condition. Let $\alpha$ be greater than
one. For a sufficiently large $a$
$Y_{a,\nu}\cap X_{\alpha a}$ is topologically a product
$T_{\nu}\times[a,\alpha a]$. For
$\sharp=abs$ or $rel$ if
$\omega\in C^{\infty}(Y_{a,\nu}\cap X_{\alpha a}, \Omega^p_X({\rho}))$ satisfies
Dirichlet condition on $T_{\nu}\times\{a\}$ and $\sharp$ on
$T_{\nu}\times\{\alpha a\}$ we will call it enjoys {\it
Dirichlet/$\sharp$-condition}.
Moreover if $\omega\in C^{\infty}(Y_{a}\cap X_{\alpha a},
\Omega^p_X({\rho}))$ satisfies Dirichlet/$\sharp$-condition for every
connected component it will be also referred that it satisfies
Dirichlet/$\sharp$-condition.\\

We will give an explicit formula of $\Delta^p$ near a cusp. First of all
notice that since Hodge Laplacian on ${\mathbb H}^3$ commutes with the
action of $\Gamma$ it preserves the decomposition:
\[
 C^{\infty}(Y_{a,\nu},\Omega^{p}_X(\rho)) =\oplus_{i=1}^{r}C^{\infty}(Y_{a,\nu},\Omega^{p}_X(\chi_{\nu,i})).
\]
Thus for a spectral problem of Hodge Laplacian on
$L^2(Y_{a,\nu},\Omega^{p}_X(\rho))$ it is sufficient to consider one on $L^2(Y_{a,\nu},\Omega^{p}_X(\chi_{\nu,i}))$.
A direct computation will show the following lemma.
\begin{lm} Let $\Delta_{T}$ be the positive Laplacian on a flat torus:
\[
 \Delta_T=-(\partial_x^2+\partial_y^2).
\]
\begin{enumerate}
\item For $f\in C^{\infty}(Y_{a,\nu},\Omega^0_X(\chi_{\nu,i}))$ we have
\[
 \Delta^0f=e^{2u}\Delta_Tf-\partial_u^2f+2\partial_uf.
\]
\item For $\omega=fdx+gdy+hdu\in
      C^{\infty}(Y_{a,\nu},\Omega^1_X(\chi_{\nu,i}))$ we have
\begin{eqnarray*}
 \Delta^1\omega &=& (e^{2u}\Delta_Tf-\partial_u^2f+2\partial_xh)dx\\
&+& (e^{2u}\Delta_Tg-\partial_u^2g+2\partial_yh)dy\\
&+& (e^{2u}\Delta_Th-\partial_u^2h+2\partial_uh-2e^{2u}(\partial_xf+\partial_yg))du.
\end{eqnarray*}
\end{enumerate}

\end{lm}
The following minimax principle will play a key role.
\begin{fact}(\cite{Reed-Simon}: The minimax principle) Let $A$ be a selfadjoint operator acting
 on a Hilbert space $H$ which is bounded below and $D(A)$ its
 domain. Then its $n$-th eigenvalue $\mu_n(A)$ is obtained by
\[
 \mu_n(A)=\inf_{\frak{M}\in {\rm Gr}_nD(A)}\sup_{0\neq v\in\frak{M}}\frac{(Av,v)}{||v||^2}.
\]
Here ${\rm Gr}_nD(A)$ is the set of $n$-dimensional subspaces of $D(A)$.
\end{fact}
Let $A$ and $B$ are selfadjoint operators bounded below which act on a Hilbert space $H$. Suppose that they have the same domain $D$ and
that $A\geq B$, i.e. $(Av,v)\geq (Bv,v)$ for any $v\in D$. Then {\bf
Fact 2.1} immediately implies

\begin{lm}
\[
 \mu_n(A)\geq \mu_n(B)
\]
\end{lm}
Let $a$ and $a^{\prime}$ be positive numbers so that $a^{\prime}\geq
a$. Extending as $0$-map on the outside $L^{2}(X_a,\Omega^p_X(\rho))$ is
embedded into $L^{2}(X_{a^{\prime}},\Omega^p_X(\rho))$ and by this
$D(\Delta^p_{dir}|_{X_a})$ is a subspace of
$D(\Delta^p_{dir}|_{X_{a^{\prime}}})$. In particular ${\rm
Gr}_n(D(\Delta^p_{dir}|_{X_a}))$ is a subset of ${\rm
Gr}_n(D(\Delta^p_{dir}|_{X_{a^{\prime}}}))$ and the minimax principle implies
\[
 \mu_{n}(\Delta_{dir}^p|_{X_{a^{\prime}}}) \leq
 \mu_{n}(\Delta_{dir}^p|_{X_{a}}).
\]
The same argument will yield the following lemma.
\begin{lm}
\begin{enumerate}
\item Let $a$ and $a^{\prime}$ be positive numbers so that $a^{\prime}\geq a$. Then we
 have
\[
 \mu_{n}(\Delta_{dir}^p|_{X_{a^{\prime}}}) \leq \mu_{n}(\Delta_{dir}^p|_{X_{a}}).
\]
\item For a positive $a$ we have
\[
 \mu_{n}(\Delta_X^p) \leq \mu_{n}(\Delta_{dir}^p|_{X_{a}}).
\]
and
\[
 \mu_{n}(\Delta_{\sharp}^p|_{X_{a}}) \leq \mu_{n}(\Delta_{dir}^p|_{X_{a}}),
\]
where $\sharp$ is $abs$ or $rel$.
\end{enumerate}
\end{lm}
Let $\Gamma_{\nu}^{*}$ be the dual lattice of $\Gamma_{\nu}$. We
will define its {\it norm} to be
\[
 ||\Gamma_{\nu}^{*}||={\rm Min}\{|\gamma|\,|\,0\neq \gamma\in\Gamma_{\nu}^{\prime}\}.
\]
Here the modulus $|\cdot|$ is taken with respect to the standard
Euclidean metric $dx^2+dy^2$.
\begin{prop} 
\[
 \mu_{1}(\Delta_{dir}^0|_{Y_{a,\nu}}) \geq e^{2a}||\Gamma_{\nu}^{*}||^2.
\]
\end{prop}
{\bf Proof.} Let us consider a nonnegative selfadjoint operator
\[
 P_a=e^{2a}\Delta_T-\partial_u^2+2\partial_u
\]
on $L^2(Y_{a,\nu},\Omega^0(\chi_{\nu,i}))$ under Dirichlet condition at
the boundary. Since 
\[
 \Delta^0-P_a=(e^{2u}-e^{2a})\Delta_T
\]
is a nonnegative operator {\bf Lemma 2.3} implies
\[
 \mu_1(\Delta^0_{dir}|_{Y_{a,\nu}})\geq \mu_1(P_a).
\]
For $f\in
C^{\infty}_c(Y_{a,\nu},\Omega^0(\chi_{\nu,i}))$ we have
\begin{eqnarray*}
\int_{Y_{a,\nu}}(P_af,f)dv_g &=&
 e^{2a}\int_{Y_{a,\nu}}\Delta_Tf\cdot\bar{f}e^{-2u}dxdydu +
 \int_{Y_{a,\nu}}|\partial_uf|^2e^{-2u}dxdydu\\
& \geq & e^{2a}\int_{Y_{a,\nu}}\Delta_Tf\cdot\bar{f}e^{-2u}dxdydu\\
&=& e^{2a}\int^{\infty}_{a}due^{-2u}\int_{T_{\nu}}\Delta_Tf\cdot\bar{f}dxdy.
\end{eqnarray*}
Let
\[
 f=\sum_{\gamma\in{\Gamma_{\nu}^{*}}}\{f_{\gamma}(u){\bf
 e}_{\gamma}(z)+f_{\gamma}^{*}(u){\bf e}_{\gamma}(\bar{z})\}, \quad {\bf
 e}_{\gamma}(z)=\exp (2\pi i\gamma z)
\]
be a Fourier expansion with respect to $T_{\nu}$-direcrion. Here notice
that by {\bf Lemma 2.1} $\gamma$ runs through nonzero
elements of $\Gamma_{\nu}^{*}$. Then
\begin{eqnarray*}
\int_{T_{\nu}}\Delta_Tf\cdot\bar{f} &=& {\rm vol}(T_{\nu})\sum_{0\neq
 \gamma \in
 \Gamma_{\nu}^{*}}|\gamma|^2\{|f_{\gamma}(u)|^2+|f_{\gamma}^{*}(u)|^2\}\\
&\geq & ||\Gamma_{\nu}^{*}||^2 {\rm vol}(T_{\nu}) \sum_{0\neq
 \gamma \in
 \Gamma_{\nu}^{*}}\{|f_{\gamma}(u)|^2+|f_{\gamma}^{*}(u)|^2\}\\
&=& ||\Gamma_{\nu}^{*}||^2\int_{T_{\nu}}|f|^2dxdy,
\end{eqnarray*}
and therefore we have obtained
\[
 \int_{Y_{a,\nu}}(P_af,f)dv_g \geq e^{2a}||\Gamma_{\nu}^{*}||^2\int_{Y_{a,\nu}}(f,f)dv_g.
\]
Now the minimax principle implies $\mu_1(P_a)\geq
e^{2a}||\Gamma_{\nu}^{*}||^2$ and the desired result has been obtained.
\begin{flushright}
$\Box$
\end{flushright}
Changing a boundary condition the above proof is still valid to prove
the following.
\begin{prop} For $\alpha>1$ and $\sharp=abs$ or $rel$, we have
\[
 \mu_{1}(\Delta_{dir\slash\sharp}^0|_{X_{\alpha a}\cap Y_{a,\nu}}) \geq e^{2a}||\Gamma_{\nu}^{*}||^2.
\]
\end{prop}
Next we will estimate $\mu_{1}(\Delta_{dir}^1|_{Y_{a,\nu}})$ from
below. First of all here are some remarks. Let us write
\[
 a=b+\beta,\quad b, \,\beta >0.
\]
Then by change of variables
\[
 u \to u+\beta
\]
$Y_{a,\nu}$ is isometric to a warped product:
\[
 \{[b,\infty)\times T_{\nu,\beta},\,du^2+e^{-2u}(dx^2+dy^2)\}
\]
where $T_{\nu,\beta}$ is a quotient of ${\mathbb
C}$ by a lattice $e^{-\beta}\Gamma_{\nu}$. In particular the dual
lattice is $e^{\beta}\Gamma_{\nu}^{*}$ and therefore if $\beta$ is
sufficiently large its norm is greater than one. Thus we may initially
assume that
\begin{equation}
 ||\Gamma_{\nu}^{*}|| >1,
\end{equation}
and by a technical reason we choose $a$ so that $e^{2a}$ is greater than
$32$. Now we will give an estimate.\\

Let $\omega=fdx+gdy+hdu$ be an element of
$C^{\infty}_c(Y_{a,\nu},\Omega^1(\chi_{\nu,i}))$. Then a computation of
{\bf Proposition 2.1} implies
\begin{equation}
\int_{Y_{a,\nu}}\Delta_Tf\cdot\bar{f}dxdydu \geq
 ||\Gamma_{\nu}^{*}||^2\int_{Y_{a,\nu}}|f|^2dxdydu \geq \int_{Y_{a,\nu}}|f|^2dxdydu,
\end{equation}
\begin{equation}
\int_{Y_{a,\nu}}\Delta_Tg\cdot\bar{g}dxdydu \geq
 ||\Gamma_{\nu}^{*}||^2\int_{Y_{a,\nu}}|g|^2dxdydu \geq \int_{Y_{a,\nu}}|g|^2dxdydu,
\end{equation}
and
\begin{equation}
\int_{Y_{a,\nu}}\Delta_Th\cdot\bar{h}e^{-2u}dxdydu \geq
 ||\Gamma_{\nu}^{*}||^2\int_{Y_{a,\nu}}|h|^2e^{-2u}dxdydu.
\end{equation}
Using the fact
\[
 ||dx||=||dy||=e^u, \quad ||du||=1
\]
and {\bf Lemma 2.2}, an integration by parts shows
\begin{eqnarray*}
\int_{Y_{a,\nu}}(\Delta^1\omega,\omega)dv_g &=&
 \int_{Y_{a,\nu}}e^{2u}(\Delta_Tf\cdot\bar{f}+\Delta_Tg\cdot\bar{g})dxdydu\\
&+ &\int_{Y_{a,\nu}}|\nabla_Th|^2dxdydu\\
&+&
 \int_{Y_{a,\nu}}(|\partial_uf|^2+|\partial_ug|^2+|\partial_uh|^2e^{-2u})dxdydu\\
&+& 2\int_{Y_{a,\nu}}\{(\partial_xh\cdot\bar{f}+\partial_x\bar{h}\cdot f)+(\partial_yh\cdot\bar{g}+\partial_y\bar{h}\cdot g)\}dxdydu
\end{eqnarray*} 
\begin{eqnarray}
&=&
 \int_{Y_{a,\nu}}(e^{2u}-16)(\Delta_Tf\cdot\bar{f}+\Delta_Tg\cdot\bar{g})dxdydu\\
&+&
 16\int_{Y_{a,\nu}}\{(\Delta_Tf\cdot\bar{f}-|f|^2)+(\Delta_Tg\cdot\bar{g}-|g|^2)\}dxdydu\\
&+&
 \frac{1}{4}\int_{Y_{a,\nu}}\{64|f|^2+8(\partial_xh\cdot\bar{f}+\partial_x\bar{h}\cdot
 f)+|\nabla_Th|^2\}dxdydu\\
&+&
 \frac{1}{4}\int_{Y_{a,\nu}}\{64|g|^2+8(\partial_yh\cdot\bar{g}+\partial_y\bar{h}\cdot
 g)+|\nabla_Th|^2\}dxdydu\\
&+& \frac{1}{2}\int_{Y_{a,\nu}}|\nabla_Th|^2dxdydu\\
&+& \int_{Y_{a,\nu}}(|\partial_uf|^2+|\partial_ug|^2+|\partial_uh|^2e^{-2u})dxdydu.
\end{eqnarray}
By (3) and (4), (7) is nonnegative and 
\[
 (8|f|-|\nabla_Th|)^2 \leq 64|f|^2+8(\partial_xh\cdot\bar{f}+\partial_x\bar{h}\cdot
 f)+|\nabla_Th|^2
\] 
and
\[
 (8|g|-|\nabla_Th|)^2 \leq 64|g|^2+8(\partial_yh\cdot\bar{g}+\partial_y\bar{h}\cdot
 g)+|\nabla_Th|^2
\]
imply (8) and (9) are also nonnegative. Since
\begin{eqnarray*}
 \int_{Y_{a,\nu}}(e^{2u}-16)(\Delta_Tf\cdot\bar{f}+\Delta_Tg\cdot\bar{g})dxdydu
  &=& \int_{Y_{a,\nu}}(e^{2u}-16)(|\nabla_Tf|^2+|\nabla_Tg|^2)dxdydu\\
& \geq& (e^{2a}-16)\int_{Y_{a,\nu}}(|\nabla_Tf|^2+|\nabla_Tg|^2)dxdydu\\
&=& (e^{2a}-16)\int_{Y_{a,\nu}}(\Delta_Tf\cdot\bar{f}+\Delta_Tg\cdot\bar{g})dxdydu
\end{eqnarray*}
we obtain
\begin{eqnarray*}
\int_{Y_{a,\nu}}(\Delta^1\omega,\omega)dv_g &\geq &
 (e^{2a}-16)\int_{Y_{a,\nu}}(\Delta_Tf\cdot\bar{f}+\Delta_Tg\cdot\bar{g})dxdydu\\
&+&  \frac{1}{2}e^{2a}\int_{Y_{a,\nu}}|\nabla_Th|^2e^{-2u}dxdydu\\
&\geq & \frac{1}{2}e^{2a}\int_{Y_{a,\nu}}(\Delta_Tf\cdot\bar{f}+\Delta_Tg\cdot\bar{g}+\Delta_Th\cdot\bar{h}e^{-2u})dxdydu.
\end{eqnarray*}
Here we have used the fact $e^{2a}$ is greater than $32$. Thus (3), (4)
and (5) implies
\[
 \int_{Y_{a,\nu}}(\Delta^1\omega,\omega)dv_g \geq \frac{1}{2}e^{2a}||\Gamma_{\nu}^{*}||^2\int_{Y_{a,\nu}}||\omega||^2dv_g
\]
and the following proposition is a direct consequence of the minimax principle.
\begin{prop} For a sufficiently large $a$ we have
\[
 \mu_{1}(\Delta_{dir}^1|_{Y_{a,\nu}}) \geq \frac{1}{2}e^{2a}||\Gamma_{\nu}^{*}||^2.
\]
\end{prop}

Changing a boundary condition the previous computation will also yield
the following.
\begin{prop}Suppose $\alpha >1$. Then for a sufficiently large $a$ and
 $\sharp=abs$ or $rel$, we have
\[
 \mu_{1}(\Delta_{dir\slash/\sharp}^1|_{X_{\alpha a}\cap Y_{a, \nu}}) \geq \frac{1}{2}e^{2a}||\Gamma_{\nu}^{*}||^2.
\]
\end{prop}

\section{Convergence of spectrum}
As we have seen in {\bf Lemma 2.4} $\mu_n(\Delta^p_{dir}|_{X_a})$ is a
monotone decreasing function of $a$, which is bounded below by
$\mu_n(\Delta^p_X)$. In this section we will show the following fact.
\begin{thm}
\[
 \lim_{a\to\infty}\mu_n(\Delta^p_{dir}|_{X_a})=\mu_n(\Delta^p_X).
\]
\end{thm}
\begin{thm} For $\sharp=abs$ or $rel$
\[
 \lim_{a\to\infty}\mu_n(\Delta^p_{\sharp}|_{X_a})=\mu_n(\Delta^p_X).
\]
\end{thm}

\begin{cor} Fot a positive $t$ we have
\[
 {\rm Tr}[e^{-t\Delta^p_X}]=\lim_{a\to\infty}{\rm Tr}[e^{-t\Delta^p_{dir}|_{X_a}}]=\lim_{a\to\infty}{\rm Tr}[e^{-t\Delta^p_{\sharp}|_{X_a}}]
\]
for $\sharp=abs$ or $rel$.
\end{cor}
Let $\chi$ be a smooth function on $X$ so that
\begin{enumerate}
\item $0 \leq \chi \leq 1.$
\item $\chi(x)=1$ on $X_a$ and vanishes on $Y_{2a}$.
\item $|\nabla \chi| \leq a^{-1}.$
\end{enumerate}
By {\bf Lemma 2.1} we know that $\Delta^p_X$ has only pure point
spectrum. Let $\varphi_i$ be its eigenform whose eigenvalue is
$\mu_i(\Delta^p_X)$ and $\frak{M}_n$ an element of ${\rm
Gr}_nD(\Delta^p_X)$ spanned by $\{\varphi_1,\cdots,\varphi_n\}$. Then
for an arbitrary $\varphi\in \frak{M}_n$ we have
\begin{equation}
\int_X ||\nabla\varphi||^2dv_g = \int_X (\Delta^{p}\varphi, \varphi)dv_g
 \leq \mu_n(\Delta^p_X)\int_X ||\varphi||^2dv_g.
\end{equation}
The LHS is
\begin{eqnarray*}
\int_X ||\nabla\varphi||^2dv_g &=& \int_X
 ||\nabla(\chi\varphi)+\nabla((1-\chi)\varphi)||^2dv_g \\
&=& \int_X||\nabla(\chi\varphi)||^2dv_g +
 \int_X||\nabla((1-\chi)\varphi)||^2dv_g\\
&+& 2{\rm Re}\int_X(\nabla(\chi\varphi), \nabla((1-\chi)\varphi))dv_g.
\end{eqnarray*}
Since
\[
 \chi(1-\chi)\leq \frac{1}{4}\quad\mbox{and}\quad |\nabla\chi|\leq \frac{1}{a}
\]
using Schwartz inequality we have
\[
 |(\nabla(\chi\varphi), \nabla((1-\chi)\varphi))| \leq (\frac{1}{a}+\frac{1}{a^2})||\varphi||^2+(\frac{1}{a}+\frac{1}{4})||\nabla\varphi||^2.
\]
Therefore (12) implies
\begin{eqnarray*}
\mu_n(\Delta^p_X)\int_X ||\varphi||^2dv_g &\geq&
 \int_X||\nabla((1-\chi)\varphi)||^2dv_g+\int_X||\nabla(\chi\varphi)||^2dv_g\\
&-&2(\frac{1}{a}+\frac{1}{a^2})\int_X||\varphi||^2dv_g
 -2(\frac{1}{a}+\frac{1}{4})\int_X||\nabla\varphi||^2dv_g\\
& \geq & \int_X||\nabla((1-\chi)\varphi)||^2dv_g\\
&-&2\{\frac{1}{a}+\frac{1}{a^2}+\mu_n(\Delta^p_X)(\frac{1}{a}+\frac{1}{4})\}\int_X||\varphi||^2dv_g.
\end{eqnarray*}
Notice that $(1-\chi)\varphi$ is contained in the domain of
$\Delta^p_{dir}|_{Y_a}$. The minimax principle and {\bf Proposition 2.1}
and {\bf Propostion 2.3} shows
\begin{eqnarray*}
\int_X||\nabla((1-\chi)\varphi)||^2dv_g &\geq&
 \mu_1(\Delta^p_{dir}|_{Y_{a}})\int_X||(1-\chi)\varphi||^2dv_g\\
&\geq& \mu_1(\Delta^p_{dir}|_{Y_{a}})\int_{Y_{2a}}||\varphi||^2dv_g\\
&\geq&Ce^{2a}\int_{Y_{2a}}||\varphi||^2dv_g,
\end{eqnarray*}
where $C$ is a positive constant independent of $a$. So we have obtained
\[
 \{\mu_n(\Delta^p_X)+2(\frac{1}{a}+\frac{1}{a^2}+\mu_n(\Delta^p_X)(\frac{1}{a}+\frac{1}{4}))\}\int_X||\varphi||^2dv_g\\
 \geq Ce^{2a}\int_{Y_{2a}}||\varphi||^2dv_g.
\]
Now putting 
\[
 \rho_n(a)=2C^{-1}e^{-2a}\{(\frac{1}{a}+\frac{1}{a^2})+\mu_n(\Delta^p_X)(\frac{1}{a}+\frac{3}{4})\}
\]
we have proved the following proposition.
\begin{prop}
For $\varphi\in\frak{M}_n$
\[
 \int_{Y_{2a}}||\varphi||^2dv_g \leq \rho_n(a)\int_X||\varphi||^2dv_g.
\]
\end{prop}
Let us fix a positive $a_0$ so that $Y_{a_0}$ is a disjoint union:
\[
 Y_{a_0}=\amalg_{\nu=1}^hT_{\nu}\times[a_0, \infty).
\]
Moreover we assume that $e^{2a_0}>32$ and that
$e^{a_0}||\Gamma_{\nu}^*||>1$ for every $\nu$. This choice guarantees
a use of {\bf Proposition 2.2} and {\bf Proposition 2.4} for an arbitrary $a$
greater than $a_0$. Let us fix such a $a$ and an any $\alpha$ greater
than two. Let $\phi_i$ be an eigenform of $\Delta^p_{\sharp}|_{X_{\alpha
a}}$ whose eigenvalue is $\mu_i(\Delta^p_{\sharp}|_{X_{\alpha
a}})$ and $\frak{M}_n(\alpha a)$ an element of ${\rm
Gr}(\Delta^p_{\sharp}|_{X_{\alpha a}})$ spanned by
$\{\phi_1,\cdots,\phi_n\}$. Then for $\phi\in \frak{M}_n(\alpha a)$ we
have
\begin{equation}
\int_{X_{\alpha a}}||\nabla\phi||^2dv_g=\int_{X_{\alpha
 a}}(\Delta^p\phi,\phi)dv_g \leq \mu_n(\Delta^p_{\sharp}|_{X_{\alpha a}})\int_{X_{\alpha a}}||\phi||^2dv_g.
\end{equation}
Using {\bf Proposition 2.2} and {\bf Proposition 2.4} instead {\bf
Proposition 2.1} and {\bf Proposition 2.3}, respectively the previous
computation will show
\begin{eqnarray*}
\mu_n(\Delta^p_{\sharp}|_{X_{\alpha a}})\int_{X_{\alpha
 a}}||\phi||^2dv_g &\geq&  \int_{X_{\alpha a}}||\nabla\phi||^2dv_g\\
 &\geq& Ce^{2a}\int_{Y_{2a}\cap
 X_{\alpha a}}||\phi||^2dv_g\\
&-& 2\{\frac{1}{a}+\frac{1}{a^2}+\mu_n(\Delta^p_{\sharp}|_{X_{\alpha
 a}})(\frac{1}{a}+\frac{1}{4})\}\int_{X_{\alpha a}}||\phi||^2dv_g,
\end{eqnarray*}
which yields
\[
 \int_{Y_{2a}\cap X_{\alpha a}}||\phi||^2dv_g \leq 2C^{-1}e^{-2a}\{(\frac{1}{a}+\frac{1}{a^2})+\mu_n(\Delta^p_{\sharp}|_{X_{\alpha
 a}})(\frac{1}{a}+\frac{3}{4})\}\int_{X_{\alpha a}}||\phi||^2dv_g.
\]
Since by {\bf Lemma 2.4} we know 
\[
 \mu_n(\Delta^p_{\sharp}|_{X_{\alpha a}}) \leq \mu_n(\Delta^p_{dir}|_{X_{a_0}})
\]
we have proved the following.
\begin{prop} Suppose $a$ is greater than $a_0$. Then for $\phi\in \frak{M}_n(\alpha a)$ we have
\[
 \int_{Y_{2a\cap X_{\alpha a}}}||\phi||^2dv_g \leq
 \rho^{0}_n(a)\int_{X_{\alpha a}}||\phi||^2dv_g,
\]
where
\[
 \rho^{0}_n(a)=2C^{-1}e^{-2a}\{(\frac{1}{a}+\frac{1}{a^2})+\mu_n(\Delta^p_{dir}|_{X_{a_0}})(\frac{1}{a}+\frac{3}{4})\}.
\]
\end{prop}
{\bf A proof of Theorem 3.1}\\

As before let $\varphi_i$ be an eigenvector of $\Delta^p_X$ whose
eigenvalue is $\mu_i(\Delta^p_X)$ and $\frak{M}_{n,\chi}$ an
$n$-dimensional subspace of $D(\Delta^p_{dir}|_{X_{2a}})$ spanned by
$\{\chi\varphi_1,\cdots,\chi\varphi_n\}$. Let us choose
$\varphi\in\frak{M}_n$ such that
\[
 \frac{\int_X||\nabla(\chi\varphi)||^2dv_g}{\int_X||\chi\varphi||^2dv_g}=\sup_{f\in\frak{M}_{n,\chi}}\frac{\int_X||\nabla
 f||^2dv_g}{\int_X||f||^2dv_g}.
\]
Since by the minimax principle the RHS is greater than or equal to
$\mu_n(\Delta^p_{dir}|_{X_{2a}})$ we know
\[
 \int_X||\nabla(\chi\varphi)||^2dv_g \geq \mu_n(\Delta^p_{dir}|_{X_{2a}})\int_X||\chi\varphi||^2dv_g.
\]
On the other hand by a choice of $\chi$ we have
\begin{eqnarray*}
||\nabla(\chi\varphi)||^2 &\leq& ||\nabla\chi\cdot\varphi||^2+2|{\rm
 Re}(\nabla\chi\cdot\varphi,\chi\nabla\varphi)|+||\chi\nabla\varphi||^2\\
&\leq& \frac{1}{a^2}||\varphi||^2+\frac{2}{a}|{\rm Re}(\varphi,\nabla\varphi)|+||\nabla\varphi||^2\\
&\leq& (\frac{1}{a^2}+\frac{1}{a})||\varphi||^2+(\frac{1}{a}+1)||\nabla\varphi||^2
\end{eqnarray*}
and therefore
\begin{eqnarray*}
& &(\frac{1}{a^2}+\frac{1}{a})\int_X||\varphi||^2dv_g
 + (\frac{1}{a}+1)\int_X||\nabla\varphi||^2dv_g\\
&\geq& \mu_n(\Delta^p_{dir}|_{X_{2a}}) \int_X||\chi\varphi||^2dv_g\\
&\geq&  \mu_n(\Delta^p_{dir}|_{X_{2a}}) \int_{X_{a}}||\varphi||^2dv_g\\
&=& \mu_n(\Delta^p_{dir}|_{X_{2a}})(\int_{X}||\varphi||^2dv_g-\int_{Y_{a}}||\varphi||^2dv_g).
\end{eqnarray*}
Using (12) {\bf Proposition 3.1} yields
\[
 (\frac{1}{a}+1)\mu_n(\Delta^p_X)\int_{X}||\varphi||^2dv_g \geq \{\mu_n(\Delta^p_{dir}|_{X_{2a}})(1-\rho_n(\frac{a}{2}))-(\frac{1}{a^2}+\frac{1}{a})\}\int_{X}||\varphi||^2dv_g
\]
and in particular
\[
 (\frac{1}{a}+1)\mu_n(\Delta^p_X)\geq \mu_n(\Delta^p_{dir}|_{X_{2a}})(1-\rho_n(\frac{a}{2}))-(\frac{1}{a^2}+\frac{1}{a}).
\]
Now notice that
\[
 \lim_{a\to\infty}\rho_n(\frac{a}{2})=0,
\]
and that by {\bf Lemma 2.4}
\[
 \mu_n(\Delta^p_X)\leq \lim_{a\to\infty}\mu_n(\Delta^p_{dir}|_{X_{2a}}),
\]
the desired result has been obtained.
\begin{flushright}
$\Box$
\end{flushright}

{\bf A proof of Theorem 3.2.}\\

 Since a proof is almost same as one of {\bf
Theorem 3.1} we will only indicate where a modification is necessary. As
before let $\phi_i$ be an eigenform of $\Delta^p_{\sharp}|_{X_{3a}}$
whose eigenvalue is $\mu_i(\Delta^p_{\sharp}|_{X_{3a}})$ and
$\frak{M}_n(3a)_{\chi}$ a $n$-dimensional subspace of
$D(\Delta^p_{dir}|_{X_{2a}})$ spanned by
$\{\chi\phi_1,\cdots,\chi\phi_n\}$. We choose $\phi\in\frak{M}_n(3a)$ so
that
\[
 \frac{\int_{X_{3a}}||\nabla(\chi\phi)||^2dv_g}{\int_{X_{3a}}||\chi\phi||^2dv_g}=\sup_{f\in\frak{M}_{n}(3a)_{\chi}}\frac{\int_X||\nabla
 f||^2dv_g}{\int_X||f||^2dv_g}.
\]
Then the minimax principle again shows
\[
  \int_{X_{3a}}||\nabla(\chi\phi)||^2dv_g \geq
  \mu_n(\Delta^p_{dir}|_{X_{2a}})\int_{X_{3a}}||\chi\phi||^2dv_g.
\]
Using (13) and {\bf Proposition 3.2} instead (12) and {\bf Proposition
3.1}, respectively the same computation as in {\bf Theorem 3.1} will yield
\[
  (\frac{1}{a}+1)\mu_n(\Delta^p_{\sharp}|_{X_{3a}})\geq \mu_n(\Delta^p_{dir}|_{X_{2a}})(1-\rho_n^0(\frac{a}{2}))-(\frac{1}{a^2}+\frac{1}{a}).
\]
By {\bf Lemma 2.4} $\mu_n(\Delta^p_{\sharp}|_{X_{3a}})$ is bounded by
$\mu_n(\Delta^p_{dir}|_{X_{3a}})$ from above and thus we
obtain by {\bf Theorem 3.1}
\[
 \lim_{a\to\infty}\mu_n(\Delta^p_{\sharp}|_{X_{3a}})=\lim_{a\to\infty}\mu_n(\Delta^p_{dir}|_{X_{3a}})=\mu_n(\Delta^p_X).
\]
\begin{flushright}
$\Box$
\end{flushright}

\section{A theorem of Cheeger-M\"{u}ller type}
By Hodge theory the $p$-th cohomology groups $H^{p}(X_a,\rho)$ and
$H^{p}(X,\rho)$ is isomorphic to ${\rm Ker}\,\Delta^{p}_{abs}|_{X_a}$
and ${\rm Ker}\,\Delta^{p}_{X}$, respectively. Here notice that
both of them have only pure point spectrum. Since for a sufficiently large
$a$ $H^{p}(X_a,\rho)$ is isomorphic to $H^{p}(X,\rho)$ by restriction, ${\rm Ker}\,\Delta^{p}_{X}$ is also isomorphic to
${\rm Ker}\,\Delta^{p}_{abs}|_{X_a}$. Let $h^{p}(X,\rho)$ be the
dimension of $H^{p}(X,\rho)$. For every $\nu$, since
$\rho|_{\Gamma_{\nu}}$ fixes only the zero vector, we will know
$H^p(T_{\nu},\rho)$ vanishes for every $p$ and $\nu$. Thus
$H^{\cdot}(X_a,\rho)$ is isomorphic to $H^{\cdot}(X_a,\partial
X_a,\,\rho)$. Moreover by Poincar\'{e} duality we have that
\[
 h^{p}(X,\rho)=h^{3-p}(X,\rho).
\]
In particular since our assumption implies that $ h^{0}(X,\rho)$
vanishes so does $h^{3}(X,\rho)$. Moreover Hodge $*$ operator induces an
isomorphism
\[
 {\rm Ker}\,\Delta^{p}_{abs}|_{X_a}\stackrel{*}\simeq {\rm Ker}\,\Delta^{3-p}_{rel}|_{X_a},
\]
and we have identities
\[
 h^{p}(X,\rho)={\rm dim}{\rm Ker}\,\Delta^{p}_{abs}|_{X_a}={\rm Ker}\,\Delta^{p}_{rel}|_{X_a}.
\]
{\it A partial spectral zeta function} of $\Delta_{\sharp}^{p}|_{X_a}$ and
$\Delta^{p}_X$ are defined to be
\[
 \zeta_{X_a,\sharp}^{(p)}(z,\rho)=\frac{1}{\Gamma(z)}\int^{\infty}_{0}\{{\rm Tr}[e^{-t\Delta_{\sharp}^{p}|_{X_a}}]-h^{p}(X,\rho)\}t^{z-1}dt,
\]
and
\[
 \zeta_{X}^{(p)}(z,\rho)=\frac{1}{\Gamma(z)}\int^{\infty}_{0}\{{\rm Tr}[e^{-t\Delta^{p}_X}]-h^{p}(X,\rho)\}t^{z-1}dt,
\]
respectively. (Here $a$ is assumed to be sufficiently large.) If ${\rm Re}\,z$ is sufficiently large they absolutely
converge and are meromorphically continued to the whole plane. Moreover
they are regular at the origin. Since
Hodge $*$ operator commutes with Laplacian,
\begin{equation}
 \zeta_{X}^{(p)}(z,\rho)=\zeta_{X}^{(3-p)}(z,\rho).
\end{equation}
and also since it interchanges two boundary conditions, we have
\begin{equation}
 \zeta_{X_a,abs}^{(p)}(z,\rho)=\zeta_{X_a,rel}^{(3-p)}(z,\rho).
\end{equation}
Now {\it a spectral zeta function} of $X_a$ and $X$ are defined to be
\[
  \zeta_{X_a}(z,\rho)=\sum_{p=0}^{3}(-1)^p p\cdot \zeta_{X_a,abs}^{(p)}(z,\rho),
\]
and
\[
   \zeta_{X}(z,\rho)=\sum_{p=0}^{3}(-1)^p p\cdot \zeta_{X}^{(p)}(z,\rho),
\]
respectively. Note that (15) and (16) imply
\[
 \zeta_{X_a}(z,\rho)=2\zeta_{X_a,rel}^{(1)}(z,\rho)-\zeta_{X_a,abs}^{(1)}(z,\rho)-3\zeta_{X_a,rel}^{(0)}(z,\rho)
\]
and
\[
 \zeta_{X}(z,\rho)=\zeta_{X}^{(1)}(z,\rho)-3\zeta_{X}^{(0)}(z,\rho).
\]
Lebesgue's convergence theorem and {\bf Corollary 3.1} yields the
following.
\begin{thm}
Suppose ${\rm Re}\,z$ is sufficiently large. Then
\[
 \lim_{a\to \infty} \zeta_{X_a}^{\prime}(z,\rho)=\zeta_{X}^{\prime}(z,\rho).
\]
\end{thm}
In this section we will show the following theorem. 
\begin{thm}
\[
 \lim_{a\to \infty}\zeta_{X_a}^{\prime}(0,\rho)=\zeta_{X}^{\prime}(0,\rho).
\]
\end{thm}
Since the origin is in the outside of the region of absolutely
convergence it will need an extra care.\\

For a finite dimensional vector space $V$ we set
\[
 {\rm det}V=\wedge^{{\rm dim} V}V
\]
and {\it the determinant} of a bounded complex of finite dimensional
vector spaces $(C^{\cdot},\,\partial)$ is defined to be
\[
 {\rm det}(C^{\cdot},\,\partial)=\otimes_{i} ({\rm det}C^{i})^{(-1)^{i}}.
\] 
Here for a complex vector space $L$ of dimension one $L^{-1}$ is its
dual. By Knudsen and Mumford it is known that there is a canonical isomorphism
\begin{equation}
 {\rm det}(C^{\cdot},\,\partial)\simeq \otimes_{i}{\rm det}\,H^{i}(C^{\cdot},\,\partial)^{(-1)^i}.
\end{equation}
Let $\Sigma=\{\Sigma_p\}_{p}$ be a triangulation of $X_a$ where
$\Sigma_p$ is the set of $p$-simplices and ${\bf e}=\{{\bf
e}_1,\cdots,{\bf e}_r\}$ a unitary base of $\rho$. We define a Hermitian
inner product on the group of $p$-cochains:
\[
 C^{p}(\Sigma,\rho)=C^{p}(\Sigma)\otimes \rho,
\]
so that $\{[\sigma]^{*}\otimes{\bf e}_i\}$ form its unitary base, where
$[\sigma]^{*}$ is the dual vector of $[\sigma]$. Now (16) induces a
metric $||\cdot||_{FR,a}$ on ${\rm
det}H^{\cdot}(X_a,\rho)=\otimes_{i}{\rm
det}\,H^{i}(C^{\cdot}(\Sigma,\rho))^{(-1)^i}$, which is {\it
Franz-Reidemeister metric} by definition. Thus by the isomorphism
\[
 H^{\cdot}(X_a,\rho)\simeq H^{\cdot}(X,\rho),
\]
we may regard ${\rm det}H^{\cdot}(X,\rho)$ a one dimensional complex
vector space with a metric $||\cdot||_{FR,a}$. Notice that they are
independent of $a$ as far as it is sufficiently large since we can
use the same triangulations to define them. Thus its limit
\[
 ||\cdot||_{FR}=\lim_{a\to \infty}||\cdot||_{FR,a}
\]
is well-defined. For a later purpose we will describe it in terms of a
combinatric zeta function. A triangulation $\Sigma$ of $X_a$ induces one $\tilde{\Sigma}$ on
the universal covering $\tilde{X_a}$ and the former may be a
quotient of the latter by the action of the fundamental group
$\Gamma$. Let $\{\sigma^{(p)}_1,\cdots,\sigma^{(p)}_{\gamma_{p}}\}$ the
set of $p$-simplices. Then $C_{p}(\tilde{\Sigma})$ is a free ${\mathbb
C}[\Gamma]$-module genereted by these elements. A twisted chain complex
is defined to be
\[
 C_{\cdot}(\Sigma,\rho)=C_{\cdot}(\tilde{\Sigma})\otimes_{{\mathbb C}[\Gamma]}\rho,
\]
which is a bounded complex of finite dimensional vector spaces. We will
introduce a Hermitian inner product so that $\{\sigma^{(p)}_i\otimes{\bf
e}_j\}$ is a unitary base. Here is an explict description of the
boundary map: Let $\sum_{k}(-1)^k\gamma_k[\sigma^{(p-1)}_{i_k}]\, (\gamma_k\in \Gamma)$ be the boundary of
$[\sigma^{(p)}_i]\in C^{p}(\tilde{\Sigma})$. Then 
\[
 \partial([\sigma^{(p)}_i]\otimes{\bf e}_j)=\sum_{k}(-1)^k
 [\sigma^{(p-1)}_{i_k}]\otimes\rho(\gamma_k){\bf e}_j.
\]
Let $(C^{\cdot}(\Sigma,\rho),\,\delta)$ be the dual complex. By the
inner product we may identify $C^{\cdot}(\Sigma,\rho)$ with
$C_{\cdot}(\Sigma,\rho)$ and in particular the dual vector of
$[\sigma^{(p)}_i]\otimes{\bf e}_j$ will be identified with itself. Thus
$(C^{\cdot}(\Sigma,\rho),\,\delta)$ is a complex such that
$C^{p}(\Sigma,\rho)$ is nothing but
$C_{p}(\Sigma,\rho)$ althogh the differential $\delta$ is the
Hermitian dual of $\partial.$ Let us define {\it a (positive) combinatric
Laplacian} $\Delta^p_{comb}$ on $C_p(\Sigma, \rho)=C^p(\Sigma, \rho)$ to be
\[
 \Delta^p_{comb}=\partial \delta+\delta \partial.
\]
Then we have
\[
 H^{p}(X_a,\rho)=H_{p}(X_a,\rho)={\rm Ker}[\Delta_{comb}^p],
\]
and both of them have the same inner product $(\cdot,\cdot)_{l^2,X_a}$ induced
by one of $C_{p}(\Sigma,\rho)$. It induces a metric
$|\cdot|_{l^2,X_a}$ on their determinant $\otimes_{p}{\rm
det}H^{p}(X_a,\rho)^{(-1)^p}$ and $\otimes_{p}{\rm det}H_{p}(X_a,\rho)^{(-1)^p}$. {\it A
combinatric zeta function} is defined as
\[
  \zeta_{comb}(s,X_a)=\sum_{p}(-1)^p p\cdot \zeta^{(p)}_{comb}(s,X_a),
\]
where
\[
 \zeta^{(p)}_{comb}(s,X_a)=\sum_{\lambda}\lambda^{-s}.
\]					       
Here $\lambda$ runs through positive eigenvalues of $\Delta^p_{comb}$ on
$C^{p}(\Sigma,\rho)$. By definition {\it a modified Franz-Reidemeister torsion}
$\tau^{*}(X_a,\rho)$ is 
\[
 \tau^{*}(X_a,\rho)=\exp (-\frac{1}{2}\zeta^{\prime}_{comb}(0,\rho)).
\]
Notice that if $H^{1}(X,\rho)$ vanishes so does every
$H^{p}(X,\rho)$ by Poincar\'{e} duality and our torsion is nothing
but the usual
Franz-Reidemeister torsion $\tau(X_a,\rho)$(\cite{Ray-Singer}). Now it is
known that $||\cdot||_{FR,a}$ is
equal to $|\cdot|_{l^2,X_a}\cdot\tau^{*}(X_a,\rho)$(\cite{Bismut-Zhang}\cite{Ray-Singer}). By construction since both
$|\cdot|_{l^2,X_a}$ and $\tau^{*}(X_a,\rho)$ depend only on a
triangulation $\Sigma$ we know they are independent of $a$ as far as it
is sufficiently large. Thus putting 
\[
 |\cdot|_{l^2,X}=\lim_{a\to \infty}|\cdot|_{l^2,X_a}, \quad
  \tau^{*}(X,\rho)=\lim_{a\to \infty}\tau^{*}(X_a,\rho),
\] 
we have
\[
 ||\cdot||_{FR}=|\cdot|_{l^2,X}\cdot\tau^{*}(X,\rho).
\]

On the other hand since $H^{p}(X_a,\rho)$ is isomorphic to 
\[
 {\rm Ker} \Delta_{abs}^{p}|_{X_a}\subset L^{2}(X_a,\Omega^{p}(\rho))
\]
the inner product on $L^{2}(X_a,\Omega^{p}(\rho))$ induces one on
$H^{p}(X_a,\rho)$. Thus by the isomorphism $H^{p}(X,\rho)\simeq
H^{p}(X_a,\rho)$ we have a metric $|\cdot|_{L^2,X_a}$ on ${\rm det}H^{\cdot}(X,\rho)$
and {\it Ray-Singer metric} $||\cdot||_{RS,a}$ is defined to be
\[
 ||\cdot||_{RS,a}=|\cdot|_{L^2,X_a}\cdot {\rm exp}(-\frac{1}{2}\zeta^{\prime}_{X_a}(0,\,\rho)).
\]
Similary using the canonical isomorphism
\[
 H^{p}(X,\rho)\simeq {\rm Ker} \Delta^{p}_X\subset L^{2}(X,\Omega^{p}(\rho)),
\]
{\it Ray-Singer metric} $||\cdot||_{RS}$ on ${\rm det}H^{\cdot}(X,\rho)$
is defined as
\[
 ||\cdot||_{RS}=|\cdot|_{L^2,X}\cdot {\rm exp}(-\frac{1}{2}\zeta^{\prime}_{X}(0,\,\rho)).
\]

\begin{prop}
\[
 \lim_{a\to\infty}|\cdot|_{L^2,X_a}=|\cdot|_{L^2,X}.
\]
\end{prop}

In fact for a sufficiently large $a$ let
$\{\phi_{a,i}\}_{i}$ be an orthonormal base of ${\rm
Ker}\,\Delta_{abs}^{p}|_{X_a}$ and we define a map
\[
 {\rm Ker}\Delta^{p}_X \stackrel{P_a}\to {\rm Ker}\Delta_{abs}^{p}|_{X_a}
\]
to be
\[
 P_a\psi=\sum_{i}\int_{X_a}(\psi,\phi_{a,i})dv_g\cdot \phi_{a,i}.
\]
Then we claim the following.
\begin{lm}
\[
 \lim_{a\to\infty}\int_{X_a}||\psi-P_a\psi||^2dv_g=0.
\]
\end{lm}
With {\bf Proposition 3.1} it will yield the following corollary which
will imply {\bf Proposition 4.1}.
\begin{cor} For $\psi\in {\rm Ker}\Delta^{p}_X$ we have
\[
 \lim_{a\to \infty}\int_{X_a}||P_a\psi||^2dv_g = \int_X||\psi||^2dv_g.
\]
\end{cor}
Here is a proof of the lemma.\\

{\bf Proof of Lemma 4.1}. For simplicity in the following arguments all positive constants
independent of $a$ will be denoted by $C$. Let $\phi_{\lambda}$ be
an eigenform of $\Delta_{abs}^{p}|_{X_a}$ whose eigenvalue is $\lambda$
satisfying 
\[
 \int_{X_a}||\phi_{\lambda}||^2dv_g=1.
\]
and we expand $\psi$ as
\[
 \psi=\sum_{\lambda}\int_{X_a}(\psi,\phi_{\lambda})dv_g\cdot \phi_{\lambda}.
\]
Since we have
\[
 \int_{X_a}||\psi-P_a\psi||^2dv_g = \sum_{\lambda > 0}|\int_{X_a}(\psi,\phi_{\lambda})dv_g|^2,
\]
it is sufficient to show that for $\phi=\phi_{\lambda}$ 
\[
 |\int_{X_a}(\psi,\phi)dv_g|\leq Ce^{-a}(\int_{X_a}||\psi||^2dv_g+C).
\]
Let us choose $\chi\in C^{\infty}_{c}(X_a)$ so that
\begin{enumerate}
\item $0\leq \chi \leq 1$.
\item $|\nabla \chi|$, $|\Delta\chi|$ are bounded by $1$.
\item $\chi\equiv 1$ on $X_{a/2}$.
\end{enumerate}
By Stokes theorem we have
\begin{eqnarray*}
\int_{X_a}(\Delta^{p}(\chi\psi),\phi)dv_g &=&
 \int_{X_a}(\chi\psi,\Delta^{p}\phi)dv_g\\
&=& \lambda \int_{X_a}\chi(\psi,\phi)dv_g
\end{eqnarray*}
Since $\Delta^{p}\psi=0$ and by the property 3 of $\chi$, the LHS becomes
\begin{eqnarray*}
\int_{X_a}(\Delta^{p}(\chi\psi),\phi)dv_g &=&
 \int_{X_a}(\Delta\chi\cdot\psi,\phi)dv_g+2\int_{X_a}(\nabla\chi\cdot\nabla\psi,\phi)dv_g\\
&=& \int_{Y_{a/2}\cap X_a}(\Delta\chi\cdot\psi,\phi)dv_g+2\int_{Y_{a/2}\cap X_a}(\nabla\chi\cdot\nabla\psi,\phi)dv_g
\end{eqnarray*}
and the property 2 of $\chi$ will imply
\begin{eqnarray*}
 |\int_{Y_{a/2}\cap X_a}(\Delta\chi\cdot\psi,\phi)dv_g|&\leq&
  \frac{1}{2}(\int_{Y_{a/2}\cap X_a}||\psi||^2dv_g+\int_{Y_{a/2}\cap
  X_a}||\phi||^2dv_g)\\
&\leq& \frac{1}{2}(\int_{Y_{a/2}}||\psi||^2dv_g+\int_{Y_{a/2}\cap
  X_a}||\phi||^2dv_g)
\end{eqnarray*}
On the other hand by {\bf Proposition 3.1} we have
\begin{eqnarray*}
\int_{Y_{a/2}}||\psi||^2dv_g &\leq& Ce^{-a} \int_{X}||\psi||^2dv_g\\
&=& Ce^{-a}(\int_{X_a}||\psi||^2dv_g+\int_{Y_a}||\psi||^2dv_g)\\
&\leq & Ce^{-a}\int_{X_a}||\psi||^2dv_g+Ce^{-a}\int_{Y_{a/2}}||\psi||^2dv_g,
\end{eqnarray*}
and therefore changing $C$ we obtain
\[
 \int_{Y_{a/2}}||\psi||^2dv_g \leq Ce^{-a}\int_{X_a}||\psi||^2dv_g.
\] 
Using {\bf Proposition 3.2} instead {\bf Proposition 3.1} the same
computation will show
\[
 \int_{Y_{a/2}\cap X_a}||\phi||^2dv_g \leq Ce^{-a}\int_{X_a}||\phi||^2dv_g=Ce^{-a}
\]
and thus
\[
 |\int_{Y_{a/2}\cap X_a}(\Delta\chi\cdot\psi,\phi)dv_g| \leq Ce^{-a}(\int_{X_a}||\psi||^2dv_g+C).
\]

Next we will estimate the second term. Using the property 2 of $\chi$ we
have
\[
 |2\int_{X_a}(\nabla\chi\cdot\nabla\psi,\phi)dv_g|\leq
  \int_{Y_{a/2}}||\nabla\psi||^2dv_g+ \int_{Y_{a/2}\cap X_a}||\phi||^2dv_g
\]
and since
\[
  \int_{Y_{a/2}}||\nabla\psi||^2dv_g \leq  \int_{X}||\nabla\psi||^2dv_g=\int_{X}(\psi,\Delta^{p}\psi)dv_g=0,
\]
it is bounded by $Ce^{-a}$. Let us consider the RHS. The property 3 implies
\[
 \int_{X_a}\chi(\psi,\phi)dv_g=\int_{X_{a/2}}(\psi,\phi)dv_g+\int_{Y_{a/2}\cap X_a}\chi(\psi,\phi)dv_g.
\]
But by the previous arguments we know the last term is bounded by
$Ce^{-a}(\int_{X_{a}}||\psi||^2dv_g+C)$. Combining all of these we will
obtain
\[
 |\int_{X_{a/2}}(\psi,\phi)dv_g| \leq Ce^{-a}(\int_{X_{a}}||\psi||^2dv_g+C).
\]
Now notice that
\begin{eqnarray*}
|\int_{X_a}(\psi,\phi)dv_g-\int_{X_{a/2}}(\psi,\phi)dv_g|&= &
 |\int_{X_a\cap Y_{a/2}}(\psi,\phi)dv_g|\\
 &\leq & \int_{X_a\cap Y_{a/2}}|(\psi,\phi)|dv_g\\
&\leq&
 \frac{1}{2}(\int_{Y_{a/2}}||\psi||^2dv_g+\int_{Y_{a/2}\cap X_a}||\phi||^2dv_g)\\
&\leq& Ce^{-a}(\int_{X_{a}}||\psi||^2dv_g+C),
\end{eqnarray*}
the desired result has been obtained since
\begin{eqnarray*}
|\int_{X_a}(\psi,\phi)dv_g| &\leq &
 |\int_{X_a}(\psi,\phi)dv_g-\int_{X_{a/2}}(\psi,\phi)dv_g|+|\int_{X_{a/2}}(\psi,\phi)dv_g|\\
&\leq & Ce^{-a}(\int_{X_{a}}||\psi||^2dv_g+C).
\end{eqnarray*}

\begin{flushright}
$\Box$
\end{flushright}
Let us choose a sufficiently large $a$ and small positive $\delta$. Let $g_{0}$ be a Riemannian metric on $X$ such that
\[
 g_0(x)=
\left\{
\begin{array}{ccc}
g(x)&{\mbox if} & x\in X_{a-\delta}\\
du^2+e^{-2a}(dx^2+dy^2)&{\mbox if} & x\in Y_{a}
\end{array}
\right.
\]
We will consider a one parameter family of metrics:
\[
 g_q=(1-q)g_0+qg, \quad 0\leq q \leq 1.
\]
Let $\{{\bf e}^0,{\bf e}^1,{\bf e}^2\}$ be an orthonormal frame of
$\Omega^1|_{X_a}$ so that ${\bf e}^0=du$. Then for $g(q)$ the second fundamental
form $h(q)$ of $\partial X_a$ and its curvature tensor $R(q)$
define an elements
\[
 \hat{h}(q)=\sum_{1\leq a,b\leq 2}h(q)_{ab}{\bf
 e}^a\otimes{\bf e}^b
\]
and
\[
 \hat{R}_0(q)=\frac{1}{4}\sum_{j,k,l}R(q)_{0jkl}{\bf e}^j\otimes({\bf
 e}^k \wedge {\bf e}^l)
\]
of $\Omega^{\cdot}|_{\partial X_a}\otimes\Omega^{\cdot}|_{\partial X_a}$,
respectively. Using Berezin integral \cite{Bismut-Zhang}, $\int^{B}$, we have an
element
\[
 \phi_a=\int^{1}_{0}dq\int^{B}\hat{h}(q)\hat{R}_0(q)\in
 \Omega^{\cdot}_{\partial X_a}.
\]
\begin{fact}(\cite{Dai-Fang})
\[
 \log
 \left(\frac{||\cdot||_{RS,a}}{||\cdot||_{FR,a}}\right)=\chi(\partial
 X_a,\rho)\log 2 +\gamma\cdot r \int_{\partial X_a}\phi_a,
\]
where $\gamma$ is an absolute constant.
\end{fact}
 Notice that the term
 $\tilde{e}(g_0,g_q)$ in the original formula vanishes because the
 dimension of $X$ is three. A direct computation will show
 that the norm of $\phi_a$ is bounded by a contant $C$ which is
 independent of $a$. Thus we obtain
\[
 |\int_{\partial X_a}\phi_a| \leq C\cdot{\rm vol}({\partial X_a}) \leq
 C^{\prime} e^{-2a},
\]
where $C^{\prime}$ is also independent of $a$. Since $\partial X_a$ is a
disjoint union of flat tori and since $\rho$ is a unitary local system
Atiyah-Singer's index theorem tells us that $\chi(\partial
X_a,\rho)$ vanishes. Thus we have proved the following proposition. 

\begin{prop}
\[
 \lim_{a\to \infty}||\cdot||_{RS,a}=||\cdot||_{FR}.
\]
\end{prop}

{\bf Proposition 4.1}, {\bf Proposition 4.2} and the definition of Ray-Singer metric will imply
that $\{\zeta^{\prime}_{X_a}(z,\rho)\}_{a}$   becomes a bounded family
of holomorphic functions on a neighborhood of the origin. Therefore
by the theorem of Ascoli-Arzela there is a subfamily
$\{\zeta^{\prime}_{X_{a_n}}(z,\rho)\}_n$ which converges to a holomorphic
function. But {\bf Theorem 4.1} shows that it should be the
restriction of $\zeta^{\prime}_{X}(z,\rho)$ and we know
\[
 \lim_{a\to\infty}\zeta^{\prime}_{X_a}(0,\rho)=\zeta^{\prime}_{X}(0,\rho).
\] 
Thus {\bf Theorem 4.2} has been proved.
The following Cheeger-M\"{u}ller
type theorem is a direct consequence of it.
\begin{thm}
$||\cdot||_{FR}$ and $||\cdot||_{RS}$ coincide. In particular
\[
 {\rm exp}(-\zeta^{\prime}_{X}(0,\rho))=\left(\frac{|\cdot|_{l^2,X}}{|\cdot|_{L^2,X}}\right)^2\tau^{*}(X,\rho)^2.
\]
\end{thm}
\section{A special value of Ruelle L-function}
Let $\Gamma_{conj}$ be the set of hyperbolic conjugacy classes of
$\Gamma$. Then there is a natural bijection between $\Gamma_{conj}$ and
the set of closed geodesics of $X$. A closed geodesic will be mentioned
as {\it prime} if it is not a positive multiple of an another one. The bijection
will determine a subset $\Gamma_{prim}$ of $\Gamma_{conj}$ which
corresponds to a subset of prime closed geodesics. For $\gamma\in
\Gamma_{conj}$ its length $l(\gamma)$ is defined to be one of
corresponding closed geodesic. Now {\it Ruelle L-function} is defined to
be
\[
 R_{X}(z,\rho)=\prod_{\gamma\in \Gamma_{prim}}{\rm det}(I_r-\rho(\gamma)e^{-zl(\gamma)})^{-1}.
\]
It is known that $R_{X}(z,\rho)$ absolutely converges if ${\rm Re}\,z$
is sufficiently large and that it is meromorphically continued to the
whole plane. Since $H^{0}(\Gamma_{\nu}, \rho)$ vanishes for every
$\nu$ by our assumption so does $H^{0}(X, \rho)$. The result of Park
will imply the following fact \cite{Park2007}:
\begin{fact}
The order of $R_{X}(z,\rho)$ at the origin is $2h^{1}(X,\rho)$ and the
 leading cofficient is
\[
 \lim_{z\to 0}z^{-2h^{1}(X,\rho)}R_{X}(z,\rho)={\rm exp}(-\zeta_{X}^{\prime}(0,\rho)).
\]
\end{fact}
Here are some remarks. In \cite{Sugiyama2} we have computed only the order
of Ruelle L-function at the origin for a unitary local system of rank
one on a hyperbolic threefold with only one cusp. Soon later J. Park
has computed the order and the leading coefficient of Ruelle L-function for an arbitrary unitary local system on an odd
dimensional complete hyperbolic manifold of finite volume. Thus {\bf
Fact 5.1} is a special case of his results.
Combining {\bf Theorem 4.3} with it we will obtain the following.
\begin{thm}
\[
 \lim_{z\to 0}z^{-2h^{1}(X,\rho)}R_{X}(z,\rho)=\left(\frac{|\cdot|_{l^2,X}}{|\cdot|_{L^2,X}}\right)^2\tau^{*}(X,\rho)^2.
\]
\end{thm}
 The ratio $|\cdot|_{l^2,X}/|\cdot|_{L^2,X}$ may be interpreted as a
 period. 
In fact let us identify $H^p(X,\rho)$ and ${\rm Ker}\Delta_X^p$ by Hodge
theory. Let $\phi^{(p)}=\{\phi^{(p)}_{1},\cdots,\phi^{(p)}_{h^p(X,\rho)}\}$ and
$\psi^{(p)}=\{\psi^{(p)}_{1},\cdots,\psi^{(p)}_{h^p(X,\rho)}\}$ be its unitary bases
with respect to $(\,,\,)_{l^2,X}$ and $(\,,\,)_{L^2,X}$,
respectively. Then using a base
$\phi_{(p)}=\{\phi_{(p),1},\cdots,\phi_{(p),h^p(X,\rho)}\}$ of $H_p(X,\rho)$
which is dual to $\phi^{(p)}$ we have an expansion
\begin{equation}
 \psi^{(p)}_{i}=\sum_{j=1}^{h^p(X,\rho)}\int_{\phi_{(p),j}}\psi^{(p)}_{i}\cdot \\\phi^{(p)}_{j}.
\end{equation}
Using these coefficients {\it a period matrix of $p$-forms} and {\it a
period of $X$} are defined to be
\[
 P(X)_p=(\int_{\phi_{(p),j}}\psi^{(p)}_{i})_{ij}
\]
and
\[
 {\rm Per}(X)=\prod_{p}|\det P(X)_p|^{(-1)^p},
\]
respectively. Then (17) implies
\[
 \psi^{(p)}_{1}\wedge \cdots \wedge \psi^{(p)}_{h^p(X,\rho)}=\det P(X)_p\cdot
 \phi^{(p)}_{1}\wedge \cdots \wedge \phi^{(p)}_{h^p(X,\rho)}.
\]
Since by definition
\[
 |\otimes_p(\psi^{(p)}_{1}\wedge \cdots \wedge \psi^{(p)}_{h^p(X,\rho)})^{(-1)^p}|_{L^2,X}=|\otimes_p(\phi^{(p)}_{1}\wedge \cdots \wedge \phi^{(p)}_{h^p(X,\rho)})^{(-1)^p}|_{l^2,X}=1,
\]
we have 
\begin{eqnarray*}
\frac{|\otimes_p(\psi^{(p)}_{1}\wedge \cdots \wedge
 \psi^{(p)}_{h^p(X,\rho)})^{(-1)^p}|_{l^2,X}}{|\otimes_p(\psi^{(p)}_{1}\wedge
 \cdots \wedge \psi^{(p)}_{h^p(X,\rho)})^{(-1)^p}|_{L^2,X}} &=& |\otimes_p(\psi^{(p)}_{1}\wedge \cdots \wedge
 \psi^{(p)}_{h^p(X,\rho)})^{(-1)^p}|_{l^2,X}\\
&=& \prod_{p}|\det P(X)_p|^{(-1)^p}\\
&=& {\rm Per}(X). 
\end{eqnarray*}
Thus {\bf Theorem 5.1} may be reformulated as follows.

\begin{thm}
\[
 \lim_{z\to 0}z^{-2h^{1}(X,\rho)}R_{X}(z,\rho)=(\tau^{*}(X,\rho)\cdot {\rm Per}(X))^2.
\]
\end{thm}
\begin{cor} Suppose that $h^{1}(X,\rho)$ vanishes. Then
\[
 R_{X}(0,\rho)=\tau(X,\rho)^{2},
\]
where $\tau(X,\rho)$ is the usual Franz-Reidemeister torsion.
\end{cor}
\section{A knot complement}
Let $K$ be a knot in $S^3$ whose complement $X_{K}$ admits a complete
hyperbolic structure of finite volume and $\rho$ a unitary local system
of rank $r$ on $X_K$. We assume that the zero is the only fixed vector
of the restriction of the associated representation of $\pi_1(X_K)$ to the fundamental group at the cusp. There is a two dimensional CW-complex $L$ to which
$X_K$ is obtained by attaching 3-cells and is a deformation retract of $X_K$. The argument of \cite{MilnorW}{\bf Lemma 7.2}
will imply the following.
\begin{lm}
\[
 \tau(X_K,\rho)=\tau(L,\rho).
\]
\end{lm}
In order to compute these terms we will clarify the chain complex
associated to $L$ and $\rho$. \\

Let
\[
 \pi_1(X_K)=<x_1,\cdots,x_n \,|\,r_1,\cdots, r_{n-1}>
\] 
be Wirtinger presentation. Here $\{x_i\}_i$ (resp. $\{r_j\}_j$) is
generators (resp. relators). We will fix a generator $t$ of
$H_1(X_K,{\mathbb Z})$ which is known to be an infinite cyclic
group. Then a group ring ${\mathbb C}[H_1(X_K,{\mathbb Z})]$ is
isomorphic to Laurent polynomial ring $\Lambda={\mathbb C}[t,t^{-1}]$
and Hurewicz map induces a ring homomorphism:
\[
 {\mathbb C}[\pi_1(X_K)]\stackrel{\epsilon}\to \Lambda.
\]
which satisfies for every $i$ 
\[
 \epsilon(x_i)=t.
\]
Also the representation $\rho$ yields a homomorphism
\[
 {\mathbb C}[\pi_1(X_K)]\stackrel{\rho}\to M_{r}({\mathbb C})
\]
and taking their tensor product we have
\[
 {\mathbb C}[\pi_1(X_K)]\stackrel{\epsilon\otimes\rho}\to M_{r}(\Lambda).
\]
Finally composing this with a homomorphism induced by the natural projection from
the free group $F_n$ of $n$-generators to $\pi_1(X_K)$ we obtain a
ring homomorphism:
\[
 {\mathbb C}[F_n]\stackrel{\Phi} \to M_{r}(\Lambda).
\]
The set of $0$-cells of $L$ consists of only one point $P_0$ and one of
1-cells is
\[
 \{x_1,\cdots,x_n\}.
\]
In order to obtain the relation it is necessary to attach 2-cells
\[
 \{y_1,\cdots,y_{n-1}\},
\]
where $y_j$ realizes the relator $r_j$. Let $\tilde{L}$ be the
universal covering of $L$ and $L_{\infty}$ an infinite cyclic covering
which corresponds ro ${\rm Ker}\epsilon$. The $p$-th chain group
$C_{p}(\tilde{L},{\mathbb C})$ is a free right ${\mathbb C}[\pi_1(X_K)]$
module generated by $P_0$ (resp. $\{x_1.\cdots.x_n\}$ or
$\{y_1,\cdots,y_{n-1}\}$) for $p=0$ (resp. $p=1$ or $p=2$) and
$C_{\cdot}(L_{\infty},\rho)$ is defined to be
\[
 C_{p}(L_{\infty},\rho)=C_{p}(\tilde{L},{\mathbb C})\otimes
 _{{\mathbb C}[{\rm Ker}\epsilon]}\rho.
\] 
Thus we have obtained a complex
\[
 C_{2}(L_{\infty},\rho)\stackrel{\partial_2}\to
 C_{1}(L_{\infty},\rho)\stackrel{\partial_1}\to C_{0}(L_{\infty},\rho),
\]
which is isomorphic to
\begin{equation}
(\Lambda^{\oplus r})^{n-1} \stackrel{\partial_2}\to (\Lambda^{\oplus
r})^{n} \stackrel{\partial_1}\to \Lambda^{\oplus r}.
\end{equation}
Using Fox free differential calculus one may compute differentials
explicitly (\cite{KL}). In fact we have 
\[
\partial_1=
\left(
\begin{array}{c}
\Phi(x_1-1)\\
\vdots\\
\Phi(x_n-1)
\end{array}
\right)
=
\left(
\begin{array}{c}
\rho(x_1)t-I_r\\
\vdots\\
\rho(x_n)t-I_r
\end{array}
\right),
\]
and 
\[
\partial_2=
\left(
\begin{array}{ccc}
\Phi(\frac{\partial r_1}{\partial x_1})& \cdots & \Phi(\frac{\partial
 r_1}{\partial x_n})\\
\vdots &\ddots & \vdots\\
\Phi(\frac{\partial r_{n-1}}{\partial x_1})& \cdots & \Phi(\frac{\partial
 r_{n-1}}{\partial x_n}).
\end{array}
\right).
\]
Here each entry is an element of $M_r(\Lambda)$. Moreover a space of chains is considered as one of row vectors and
differentials act from the right. 
It is known that a determinant of a certain entry of $\partial_1$ is not
zero(\cite{Wada}). Therefore rearranging numbers we may assume that
$\det(\rho(x_n)t-I_r)$ is not zero, which will be denoted by
$\Delta_0(t)$. Now we put
\[
 \Delta_{1}(t)={\rm det}
\left(
\begin{array}{ccc}
\Phi(\frac{\partial r_1}{\partial x_1})& \cdots & \Phi(\frac{\partial
 r_1}{\partial x_{n-1}})\\
\vdots &\ddots & \vdots\\
\Phi(\frac{\partial r_{n-1}}{\partial x_1})& \cdots & \Phi(\frac{\partial
 r_{n-1}}{\partial x_{n-1}}).
\end{array}
\right)
\]
and a ratio
\[
 \Delta_{K,\rho}(t)=\frac{\Delta_1(t)}{\Delta_0(t)}
\]
is nothing but {\it a twisted Alexander function}
(\cite{KL}\cite{Kitano}\cite{Wada}). In the following we will assume
$\Delta_1(t)$ is not zero. Since $C_{\cdot}(L,\rho)$, which is
quasi-isomorphic to $C_{\cdot}(X_K,\rho)$, is obtained by
modding out (18) by an ideal generated $(t-1)$ we have a long exact
sequence:
\begin{equation}
 \begin{array}{cccccccccc}
0&\to & H_2(L_{\infty},\rho)&\stackrel{\tau_2-id}\to & H_2(L_{\infty},\rho) &\to &
 H_2(X_K,\rho) & & \\
& \to & H_1(L_{\infty},\rho)&\stackrel{\tau_1-id}\to & H_1(L_{\infty},\rho) &\to &
 H_1(X_K,\rho) & & \\
& \to & H_0(L_{\infty},\rho) & \stackrel{\tau_0-id}\to &  H_0(L_{\infty},\rho) &\to & H_0(X_K,\rho) &\to & 0,
\end{array}
\end{equation}
where $\tau_i$ is the representation matrix of the action of $t$ on
corresponding spaces. Here notice that $\Delta_i(t)\neq 0$ implies if
tensored with ${\mathbb C}(t)$ (18) becomes acyclic. Therefore every $H_{\cdot}(L_{\infty},\rho)$ is a torsion
$\Lambda$-module and in particular they are finite dimensional vector
spaces over ${\mathbb C}$. Since $h^0(X_K,\rho)=0$ we know by the
universal coefficient theorem that 
$H_0(X_K,\rho)$ vanishes. Thus $\tau_0-id$ is an isomorphism and the
vanishing of $h^1(X_K,\rho)$ is equivalent to the fact that $\tau_1-id$ is
isomorphic.   Since $\Delta_i(t)$ differs from the characteristic
polynomial of $\tau_i$ by a unit of $\Lambda$, using (18), one will
easily see that $\Delta_1(1)\neq 0$ induces
$h^1(X_K,\rho)=h^2(X_K,\rho)=0$. Conversely (19) will also show that
$h^1(X_K,\rho)=0$ yields $\Delta_1(1)\neq 0$ and
$h^2(X_K,\rho)=0$. In \cite{Sugiyama} we have proved that the
vanishing of $h^{i}(X_K,\rho)$ for all $i$ implies
\begin{equation}
 \tau(X_K,\rho)=|\Delta_{K,\rho}(1)|.
\end{equation}
Thus {\bf Corollary 5.1} and (20) implies the following.
\begin{thm}
Suppose $h^{1}(X_K,\rho)=0$. Then 
\[
 R_{X_K}(0,\rho)=|\Delta_{K,\rho}(1)|^2.
\]
\end{thm}

Here is an example of $\rho$ such that a special value of Ruelle
L-function at the origin can be computed explicitly. Let $\xi$ be a
complex number of modulus one and 
\[
 H_{1}(X_K,{\mathbb Z})\stackrel{\rho}\to U(1)
\]
a representation defined to be
\[
 \rho(t)=\xi.
\]
Composing it with Hurewicz map
we obtain a unitary character
\[
 \pi_1(X_K)\stackrel{\rho}\to U(1),
\]
which yields a unitary local system of rank one on $X_K$. 
Since $t$ represents a meridian of the boundary of a tubular
neighborhood of $K$, if $\xi\neq 1$, the required assumption of the
representation at the cusp is satisfied. Moreover
it is known
(\cite{KL}\S 3.3):
\[
 \Delta_0(t)=1-\xi t, \quad \Delta_1(t)=A_K(\xi t),
\]
where $A_{K}(t)$ is Alexander polynomial. Now let us choose $\xi$
so that $\xi\neq 1$ and that $A_{K}(\xi)\neq 0$. Then the previous
argument shows that
$h^{1}(X_K,\rho)=0$ and by {\bf Theorem 6.1} we have the following.
\begin{cor}
\[
 R_{X_K}(0,\rho)=\left|\frac{A_K(\xi)}{1-\xi}\right|^2.
\]
\end{cor}


\vspace{10mm}
\begin{flushright}
Address : Department of Mathematics and Informatics\\
Faculty of Science\\
Chiba University\\
1-33 Yayoi-cho Inage-ku\\
Chiba 263-8522, Japan \\
e-mail address : sugiyama@math.s.chiba-u.ac.jp
\end{flushright}

\end{document}